\documentclass[english,11pt]{article}
\usepackage[T1]{fontenc}
\usepackage{babel}
\usepackage{amsmath,amsfonts,amsthm}
\usepackage{color}

\DeclareMathSymbol{\leqslant}{\mathalpha}{AMSa}{"36} 
\DeclareMathSymbol{\geqslant}{\mathalpha}{AMSa}{"3E} 
\renewcommand{\leq}{\;\leqslant\;}                   
\renewcommand{\geq}{\;\geqslant\;}                   

\newtheorem{Th}{Theorem}
\newtheorem{Le}[Th]{Lemma}
\newtheorem{Pro}[Th]{Proposition}

\newcommand{\cA}{\ensuremath{\mathcal A}}
\newcommand{\cB}{\ensuremath{\mathcal B}}

\newcommand{\cD}{\ensuremath{\mathcal D}}
\newcommand{\cE}{\ensuremath{\mathcal E}}

\newcommand{\cH}{\ensuremath{\mathcal H}}

\newcommand{\cP}{\ensuremath{\mathcal P}}

\newcommand{\cR}{\ensuremath{\mathcal R}}
\newcommand{\cS}{\ensuremath{\mathcal S}}

\newcommand{\cX}{\ensuremath{\mathcal X}}

\newcommand{\bbD}{{\ensuremath{\mathbb D}} }
\newcommand{\bbE}{{\ensuremath{\mathbb E}} }

\newcommand{\bbN}{{\ensuremath{\mathbb N}} }

\newcommand{\bbP}{{\ensuremath{\mathbb P}} }

\newcommand{\bbR}{{\ensuremath{\mathbb R}} }

\newcommand{\Om}{\Omega}
\newcommand{\E}{\bbE}

\newcommand{\N}{\bbN}
\newcommand{\Ne}{\bbN^{\ast}}
\newcommand{\R}{\bbR}
\newcommand{\bbd}{\mathbf{d}}

\newcommand{\tN}{\tilde{N}}
\newcommand{\hOm}{\hat{\Om}}
\newcommand{\hcA}{\hat{\cA}}
\newcommand{\hbbP}{\hat{\bbP}}
\newcommand{\hE}{\hat{\E}}
\newcommand{\hw}{\hat{w}}
\newcommand{\crea}{\varepsilon^+}
\newcommand{\anni}{\varepsilon^-}
\newcommand{\sbbD}{\underline{\bbD}}
\newcommand{\1}{\mathbf{1}}

\newcommand{\xs}{X_{s^-}}
\newcommand{\alu}{^{(\alpha ,u)}}

\newcommand{\bK}{\bar{K}}

\newcommand{\LP}{\cH_{\cP}}
\newcommand{\LD}{{\cH}_{\bbD}}
\newcommand{\LDP}{{\cH}_{\bbD,\cP}}
\newcommand{\LDdP}{{\cH}_{\bbD\otimes{\bbd},\cP}}

\title{Application of the lent particle method\\ to Poisson driven SDE's}\date{}
\author{Nicolas BOULEAU and Laurent DENIS}
\begin{document}
\maketitle

\date{}

\begin{abstract}
We apply the Dirichlet forms version of Malliavin calculus to stochastic differential equations with jumps.  As in the continuous case this weakens significantly the assumptions on the coefficients of the SDE. In spite of the use of the Dirichlet forms theory, this approach brings also an important simplification which was not available nor visible previously : an explicit formula giving the carr\'e du champ matrix, i.e. the Malliavin matrix. Following this formula a new procedure appears, called the lent particle method which shortens the computations both theoretically and in concrete examples. In this paper which uses the construction done in \cite{bouleau-denis} we restrict ourselves to the  existence of densities, smoothness will be studied separately.

\end{abstract}

\section{Introduction}
{
During the last twenty years a significant development of the theory of Dirichlet forms occured in the direction of improving regularity results or lightening hypotheses of Malliavin calculus \cite{malliavin}. With respect to the Malliavin analysis on Wiener space, what brings the Dirichlet forms approach is threefold: a) The arguments hold under only Lipschitz hypotheses, e.g. for regularity of solutions of stochastic differential equations cf \cite{bouleau-hirsch2}, this is due to the celebrated property that contractions operate on Dirichlet forms and the \'Emile Picard iterated scheme may be performed under the Dirichlet norm. b) A general criterion exists, the energy image density property (EID), proved on the Wiener space for the Ornstein-Uhlenbeck form, and in several other cases (but still a conjecture in general since 1986 cf \cite{bouleau-hirsch1}), which provides an efficient tool for obtaining existence of densities in stochastic calculus. c) Dirichlet forms are easy to construct in the infinite dimensional frameworks encountered in probability theory (cf \cite{bouleau-hirsch2} Chap.V) and this yields a theory of errors propagation, especially for finance and physics cf \cite{bouleau3}, but also for numerical analysis of PDE's and SPDE's cf \cite{scotti}.

Extensions of Malliavin calculus to the case of stochastic differential equations with jumps have been soon proposed and gave rise to an extensive literature. The approach is either dealing with local operators acting on the size of the jumps (cf \cite{bichteler-gravereaux-jacod} \cite{coquio} \cite{ma-rockner2} etc.) or acting on the instants of the jumps (cf \cite{carlen-pardoux} \cite{denis}) or based on the Fock space representation of the Poisson space and finite difference operators (cf \cite{nualart-vives}  \cite{picard} \cite{ishikawa-kunita} etc.). In all cases  the arguments  are somewhat intricate.

We have obtained results which simplify highly the approach with local operators cf \cite{bouleau-denis}. Based on Dirichlet forms on the general Poisson space they gather the advantages of Dirichlet forms methods and simplicity of use. It may be summarized in the following way: in order to calculate the Malliavin matrix, we add a particle to the system, compute the gradient of the functional on this particle, and take back the particle before integrating by the Poisson measure. The main formula is

\[\Gamma [F] =\int_X \anni(\gamma [\crea F ])\, dN.\]
where $\gamma$ is the carr\'e du champ operator on the state space, $\Gamma$ the carr\'e du champ operator on the Poisson space and $\crea$, $\anni$ are the operations of adding and cancelling a particle.

Let us present this method, called {\it the lent particle method}, on a simple example.
Let $Y$ be a real L\'evy process with absolutely continuous L\'evy measure $\nu =k dx$, thus such that $1+\Delta Y_s\neq 0$ almost surely. We equip the space $\R \setminus \{ 0\}$ with the Dirichlet form on $L^2 (\nu )$ with carr\'e du champ operator
\[ \gamma [u] (x)=x^2 {u'}^2 (x) {\bf 1}_{\{ |x|< \frac12\}}.\]

The link of $Y$ with the associated  random Poisson measure $N$ with intensity $dt\times\sigma$ is $\forall h\in L^2_{loc}(\mathbb{R}^+)$, $\int_0^th(s)\;dY_s=\int{\bf 1}_{[0,t]}(s)h(s)x\tilde{N}(dsdx)$ where $\tilde{N}$ is the compensated measure $N-dt\times\sigma$.

 Let us consider the Dol\'eans-Dade exponential
\begin{equation}\label{exp}\mathcal{E}(Y)_t=e^{Y_t}\prod_{s\leq t}(1+\Delta Y_s)e^{\Delta Y_s}.\end{equation}

In order to study the regularity of the pair  $(Y_t ,\mathcal{E}(Y)_t )$, we proceed as follows:
\begin{enumerate}
\item We add a particle $(\alpha ,y)$ i.e. a jump to $Y$ at time $\alpha\leq t$ with size $y$ what gives

$$\varepsilon^+_{(\alpha,y)}(\mathcal{E}(Y)_t)=e^{Y_t+y}\prod_{s\leq t}(1+\Delta Y_s)e^{\Delta Y_s}(1+y)e^{-y}=\mathcal{E}(Y)_t(1+y).$$

\item We compute $\gamma [\varepsilon^+ \mathcal{E}(Y)_t]=(\mathcal{E}(Y)_t)^2 y^2 {\bf 1}_{\{ |y|< \frac12\}}$.
\item We take back the particle before integrating w.r.t. $N$  :
\[ \varepsilon^- \gamma [\varepsilon^+ \mathcal{E}(Y)_t]=\left( \mathcal{E}(Y)_t (1+y)^{-1}\right)^2 y^2 {\bf 1}_{\{ |y|< \frac12\}}\]
and
\[ \Gamma [\mathcal{E}(Y)_t]=\sum_{\alpha \leq t}\left( \mathcal{E}(Y)_t (1+\Delta Y_{\alpha})^{-1}\right)^2 \Delta Y_{\alpha}^2 {\bf 1}_{\{ |\Delta Y_{\alpha}|<\frac12\}}.\]
\end{enumerate}
As easily seen, by a similar computation, the matrix $\Gamma$ of the pair $(Y_t ,\mathcal{E}(Y_t ))$ is given by
\[ \Gamma=\sum_{\alpha \leq t}\left(
                                \begin{array}{cc}
                                  1 &  \mathcal{E}(Y)_t (1+\Delta Y_{\alpha})^{-1}\\
                                  \mathcal{E}(Y)_t (1+\Delta Y_{\alpha})^{-1} & \left( \mathcal{E}(Y)_t (1+\Delta Y_{\alpha})^{-1}\right)^2 \\
                                \end{array}
                              \right)\Delta Y_{\alpha}^2 {\bf 1}_{\{ |\Delta Y_{\alpha}|<\frac12\}}.
\]
Hence the density of the pair $(Y_t ,\mathcal{E}(Y_t ))$ under hypotheses yielding (EID) lies on the condition
$$\mbox{dim  } \mathcal{L}\left(\left(\begin{array}{c}
1\\
\mathcal{E}(Y)_t(1+\Delta Y_\alpha)^{-1}
\end{array}\right)\quad \alpha\in JT\right)=2
$$
where $JT$ denotes the jump times of $Y$ with size less that $\frac12$ between 0 and $t$. This is guaranteed if $Y$ has an infinite L\'evy measure.

\noindent Hence, summarizing our hypotheses (see Lemma \ref{HypBC} below)

{\it Let $Y$ be a real L\'evy process whose L\'evy measure is infinite, absolutely continuous w.r.t. Lebesgue measure and whose density dominates a positive continuous function near $0$, then the pair $(Y_t,\mathcal{E}(Y)_t)$ possesses a density on $\mathbb{R}^2$.}

The aim of the present article is to apply the {\it lent particle method} to stochastic differential equations with jumps.

We shall begin by explaining the method (Sections 2.1 and 2.2). The only difference with our treatment in \cite{bouleau-denis} is that we work here on a product probability space in order to be able to put also a Brownian motion or other semi-martingales in the studied SDE. With respect to the note \cite{bouleau4} we have introduced a clearer new notation, as in \cite{bouleau-denis} the operator $\varepsilon^-$ is shared from the integration by $N$.

The SDE we consider is described in Section 2.3 with the assumptions done.

The argument dealing to the Malliavin matrix is concentrated in Section 2.6 thanks to the lent particle method.

We give some examples in Section 3 of SDE driven by L\'evy processes when the L\'evy measure is underdimensioned and the diffusion matrix is degenerated. This improves known results on  the L\'evy area (Section 3.2). We end by lightening assumptions on existence of density for a McKean-Vlasov type non-linear SDE (Section 3.3) and for stable-like processes (Section 3.4).

\section{Notations and hypotheses.}
\subsection{Dirichlet structures and Poisson measures.}
Let $(X,\cX ,\nu,\bbd,\gamma)$ be a local symmetric Dirichlet
structure which admits a carr\'e du champ operator i.e. $(X,\cX ,\nu
)$ is a measured space, $\nu$ is $\sigma$-finite and the bilinear
form
\[ e [f,g]=\frac12\int\gamma [f,g]\, d\nu,\]
is a local  Dirichlet form with domain $\bbd\subset L^2 (\nu )$ and
carr\'e du champ operator $\gamma$ (see Bouleau-Hirsch
\cite{bouleau-hirsch2}, Chap. I). We assume that for all $x\in X$,
$\{ x\}$ belongs to $\cX$ and that $\nu$ is diffuse
($\nu(\{x\})=0\;\forall x$). The structure  $(X,\cX ,\nu,\bbd,\gamma)$ is called the {\it bottom structure}.

  We are given $N$ a Poisson random measure on $[0,+\infty[
\times X$ with intensity $dt\times \nu(du)$ defined on the
probability space $(\Om_1 ,\cA_1 ,\bbP_1)$ where $\Om_1$ is the
configuration space, $\cA_1$ the $\sigma$-field generated by $N$ and
$\bbP_1$ the law of $N$. We set
$\tN =N-dt\times\nu$.

 We suppose that the bottom structure and $N$ satisfy the following hypothesis :\\

\underline{Hypothesis (H0)}. The structure $(X,\cX ,\nu,\bbd,\gamma)$ with generator $(a,\mathcal{D}(a))$ is such that there exists a subspace $H$ of $\cD (a)\bigcap L^1 (\nu )$,  dense in $L^1 (\nu )\cap L^2 (\nu )$ and such that $\forall f\in H,\;\gamma[f]\in L^2(\nu)$.\\

Hypothesis (H0) implies what we call the {\it bottom core hypothesis} in \cite{bouleau-denis} and denote (BC). It is a technical condition due to the fact that the carr\'e du champ takes its values in $L^1$ and we need a set of test functions for which it has its values in $L^2$. This condition is not so restrictive. In the case of a Poisson measure induced by a real L\'evy process whose L\'evy measure is absolutely continuous w.r.t. the Lebesgue measure, the following Lemma gives a way to fulfill it:
\begin{Le}{\label{HypBC}} Let $r\in \Ne$, $(X,\cX )=(\R^r ,\mathcal{B}(\R^r))$ and $\nu =k dx$ where $k$ is non-negative and Borelian. We are given  $\xi=(\xi_{ij})_{1\leq i,j\leq r}$ an $\R^{r\times r}$-valued and
symmetric Borel function. We assume that there exist an open set $O\subset \R^r$ and a  function $\psi$ continuous on $O$ and null on $\R^r \setminus O$ such that
\begin{enumerate}
\item $\nu (\partial O)=0$,
\item  $k>0$ on $O$ $\nu$-a.e. and is locally bounded on $O$
\item  $\xi$ is locally bounded and  locally elliptic on $O$ in the sense that for any compact subset
$K$ in $O$, there exists  positive constants $c_K$ and $C_K$ such that
\[\forall x=(x_1 ,\cdots ,x_r )\in K,\ \sum_{i,j}|\xi_{i,j} (x)|\leq C_K\makebox{ and }\sum_{i,j=1}^r \xi_{ij} (x)x_i x_j \geq
c_K |x|^2  ,\]
\item $k\geq \psi>0$ $\nu$-a.e. on $O$
\item for all $i,j\in \{ 1,\cdots ,r\}$, $\xi_{i,j}\psi$ belongs to $H^1_{loc} (O)$.
\end{enumerate}
We denote by $H$ the subspace of functions $f\in L^2 (\nu )\cap L^1 (\nu )$ such that the restriction of $f$ to $O$ belongs to $C_0^{\infty} (O )$ (i.e. $\mathcal{C}^\infty$ with compact support in $O$).
Then, the bilinear form defined by
\[\forall f,g\in H,\ e(f,g)=\sum_{i,j=1}^r \int_O \xi_{i,j}(x)\partial_i f(x)\partial_j g (x)\psi(x)\, dx\]
is closable in $L^2 (\nu)$. Its closure, $(\bbd ,e)$, is a local Dirichlet form on $L^2 (\nu)$ which admits a carr\'e du champ $\gamma$. Moreover, it satisfies
hypothesis {\rm(H0)} and property {\rm(EID)} i.e. for any $d$ and for any $\mathbb{R}^d$-valued function $U$
whose components are in the domain of the form
$$ U_*[({\det}\gamma[U,U^t])\cdot
\nu ]\ll \lambda^d $$ where {\rm det} denotes the determinant
and $\lambda^d$ the Lebesgue measure on $(\R^d ,\cB (\R^d ))$.\\
\end{Le}
\begin{proof} First of all, since $\nu (\partial O)=0$, it is clear that $H$ is dense in $L^2 (\nu )\cap L^1 (\nu )$. We have
\[ \forall f,g\in H, \ e(f,g)=-\int_{\R^r} a(f)(x)g(x)\, \nu (dx),\]
where $a(f)={\bf 1}_O (x)({k(x)})^{-1}\sum_{i,j} \partial_j (\xi_{i,j} \psi\partial_i f)$. Since $k\geq \psi$ on $O$, it is clear that for any
compact set $K\subset O$, $k^{-1}$ is bounded on $K$. From this, it is clear that $a(f)$ belongs to $L^2 (\nu )$. So, $a$ is an operator with dense domain, symmetric  and $-a$ is non-negative hence it is closable. As a consequence, the bilinear form $(H,e)$ is also closable and it is clear that its closure is a Dirichlet form
on $L^2 (\nu )$ and that the carr\'e du champ operator $\gamma$ is given by
\[ \forall f\in H,\  \gamma (f)(x)=\sum_{i,j=1}^r \xi_{i,j}(x)\partial_i f(x)\partial_j f(x)\frac{\psi(x)}{k(x)},\]
with the convention $\frac00=0$.\\
This ensures that hypothesis (H0) is satisfied.\\
The (EID) property can be proved using standard arguments (see Theorem 2 in \cite{bouleau-denis}).
\end{proof}
Hypothesis (H0) is satisfied in several other cases not so simple to describe (cf \cite{fukushima-oshima-takeda} Section 3.1). \\

When the aim is to obtain density for $\mathbb{R}^d$-valued random variables, $d>1$, {\it and only in that case}, we have to suppose additional conditions : \\
{\underline{Hypothesis (H1)}}: The structure $(X,\cX ,\nu,\bbd,\gamma)$ satisfies (EID).\\
{\underline{Hypothesis (H2)}}: $X$ admits a partition of the form:
$X=B\bigcup (\bigcup_{k=1}^{+\infty}A_k )$ where for all $k$,
$A_k\in\cX$ with $\nu (A_k )<+\infty$ and $\nu (B)=0$, in such a way
that for any $k\in\Ne$ may be defined a local  Dirichlet structure
with carr\'e du champ:
\[ \cS_k =(A_k ,\cX_{|A_k}, \nu_{|A_k }, \bbd_k,\gamma_k),\]
with
$\forall f\in\bbd, \, f_{|A_k }\in\bbd_k \makebox{ and } \gamma
[f]_{|A_k}=\gamma_f [f_{|A_k}].$\\
 {\underline{Hypothesis (H3)}}: Any finite product of structures
 $\cS_k$ satisfies (EID).\\

 The need for hypotheses (H1) to (H3) lies in the argument followed in \cite{bouleau-denis} to prove (EID) on the Poisson space. As mentionned above, no case is known where (EID) fails for a local Dirichlet structure with carr\'e du champ. In the classical applications hypotheses (H0) to (H3) are fulfilled.\\

We consider also another  probability space $(\Om_2,\cA_2 ,\bbP_2)$ on which an $\R^n$-valued semimartingale $Z=(Z^1 ,\cdots,Z^n)$ is defined, $n\in\Ne$.
We adopt the following assumption  on the bracket of $Z$ and on the total variation of its finite variation part. It is satisfied if both are dominated by the Lebesgue measure, uniformly w.r.t. $w$:\\
 {\underline{Assumption on $Z$}}:

There exists a positive constant $C$ such that  for any square integrable $\R^n$-valued predictable process $h$:
\begin{equation}\label{CondCrochet} \forall t\geq 0 ,\ \E  [(\int_0^t h_s dZ_s )^2]\leq C^2 \E [\int_0^t |h_s|^2 ds].\end{equation}

The presence of this semimartingale and of the following  framework is due to the form of the SDE that we will study below cf equation (\ref{eq}). This form will be discussed in \S 2.3.

We shall work on the product probability space:
\[(\Om,\cA ,\bbP)=(\Om_1 \times \Om_2 ,\cA_1 \otimes \cA_2 ,\bbP_1
\times \bbP_2 ).\]
As in \cite{bouleau-denis}, starting from the
 Dirichlet structure on the {\it bottom space} $(X,\cX ,\nu )$ we
construct a Dirichlet form on $L^2 (\Om_1 ,\cA_1,\bbP_1 )$ and then by
considering the product of this Dirichlet structure with the trivial
one on $L^2 (\Om_2 ,\cA_2 ,\bbP_2 )$, we obtain a Dirichlet
structure $(\Om ,\cA ,\bbP ,\bbD ,\cE ,\Gamma)$ with domain
$\bbD\subset L^2 (\Om ,\cA ,\bbP)$ and carr\'e du champ operator
$\Gamma$. As we assume (H0) to (H3) we know that it satisfies
(EID) (cf \cite{bouleau-denis} or Song \cite{song}).\\
Finally we denote by $(\cA_t )_{\geq 0}$ the natural filtration of
the Poisson random measure $N$ on $[0,+\infty[\times X.$
 \subsection{Expression for the gradient and the carr\'e du
champ operator} Following \cite{bouleau-denis} \S 3.2.2, we are given an
auxiliary probability space $(R,\cR ,\rho)$ such that the dimension
of the vector space $L^2 (R,\cR
 ,\rho)$ is infinite and we  construct a
random Poisson measure $N\odot\rho$ on $[0,+\infty[\times X\times R
$ with compensator $dt\times\nu\times \rho$ such that if $N=\sum_i
\varepsilon_{(\alpha_i ,u_i )}$ then $N\odot\rho =\sum_i
\varepsilon_{(\alpha_i ,u_i ,r_i )}$ where $(r_i )$ is a sequence of
i.i.d. random variables independent of $N$ whose common law is
$\rho$ and defined on some probability space $(\hOm, \hcA ,\hbbP)$
so that $N\odot\rho$ is defined on the product probability space:
$(\Om,
\cA ,\bbP)\times (\hOm, \hcA ,\hbbP)$.\\
We assume that the Hilbert space $\bbd$ is separable so that the bottom
Dirichlet structure admits a gradient operator, $D$, and we choose a
version of it with values in
 the  space $L_0^2 (R,\cR ,\rho)=\{ g\in L^2 (R,\cR ,\rho); \int_R g(r)\rho(dr)=0\}$.\\
  We denote it by $\flat$.\\
  Let us recall some important properties:
  \begin{itemize}
  \item $\forall u\in\bbd,\ Du=u^{\flat} \in L^2 (X\times R ,\cX \otimes \cR
 ,\nu\times\rho).$
 \item $\forall u\in \bbd$, $\int_R \| u^{\flat}\|^2 (\cdot ,r)\rho (dr) =\gamma[u]$.
 \item (chain rule in dim 1) if
$F:\R\rightarrow \R$ is Lipschitz  then
$\forall u\in\bbd,\ (F\circ u)^{\flat}=(F'\circ u )u^{\flat}.$
\item (chain rule in dim $d$) if $F$ is $\mathcal{C}^1$ (continuously differentiable) and Lipschitz from $\R^d$ into
$\R$ then
\[ \forall u=(u_1 ,\cdots ,u_d) \in \bbd^d ,\ (F\circ
u)^{\flat}=\sum_{i=1}^d (F'_i \circ u ) u_i^{\flat}.\]
\end{itemize}
Finally,  although not necessary, we assume for simplicity that
constants belong to $\bbd_{loc}$ (see Bouleau-Hirsch
\cite{bouleau-hirsch2} Chap. I Definition 7.1.3.) and that
\begin{equation}\label{232}
1\in \bbd_{loc} \makebox{ which implies }\ \gamma [1]=0 \makebox{
and  } 1^{\flat}=0.
\end{equation}
We now introduce the creation and annihilation operators $\crea$ and
$\anni$:
\[\begin{array}{l} \forall (t,u)\in [0,+\infty[\times X,\forall w_1\in\Om_1,\
\crea_{(t ,u)} (w_1)=w_1{\bf 1}_{\{ (t ,u)\in supp
\, w_1\}}+(w_1+\varepsilon_{(t ,u)}\}) {\bf 1}_{\{ (t ,u)\notin supp \, w_1\}}\\
\forall (t,u)\in [0,+\infty[\times X,\forall w_1\in\Om_1,\ \anni_{(t
,u)} (w_1)=w_1{\bf 1}_{\{ (t ,u)\notin supp \,
w_1\}}+(w_1-\varepsilon_{(t ,u)}\}) {\bf 1}_{\{ (t ,u)\in supp \,
w_1\}}.
\end{array}\]
In a natural way, we extend these operators on $\Omega$ by setting
if $w=(w_1 ,w_2)$:
\[ \crea_{(t ,u)}(w)=(\crea_{(t ,u)}(w_1),w_2)\
\makebox{and}\ \anni_{(t ,u)}(w)=(\anni_{(t ,u)}(w_1),w_2),\] and
then to the functionals by
\[ \crea H(w,t ,u)=H(\crea_{(t ,u)} w,t, u)\quad\makebox{ and }\quad\anni H(w,t ,u)=
H(\anni_{(t , u)}w,t,u).\] We now recall the main Theorem of
\cite{bouleau-denis} which gives an explicit formula for a
gradient of the upper structure $(\Om ,\cA ,\bbP ,\bbD ,\cE ,\Gamma)$.\\
Let us introduce some notations. The space $\cD_0$ was introduced in hypothesis (H0). We denote $\sbbD$ the
completion of $\cD_0\otimes L^2 ([0,+\infty[,dt)\otimes \bbd$ with
respect to the norm
\begin{eqnarray*}\| H\|_{\sbbD}\!\!&=&\!\!\left( \E\int_0^{\infty}\!\!\!\int_{ X }\anni(\gamma
[H])(w,t,u)N(dt,du)\right)^{\frac12} \!+\E\int_0^{\infty}\!\!\!\int_X
(\anni|H|)(w,t,u)\eta (t,u)N(dt,du)\\
\!\!&=&\!\! {\left( \E\int_0^{\infty}\int_{ X }\gamma
[H](w,t,u)\nu (du)dt\right)^{\frac12}\! +\E\int_0^{\infty}\int_X
|H|(w,t,u)\eta (t,u)\nu(du)dt},\end{eqnarray*} where $\eta$ is a fixed
positive function in $L^2 (\R^+ \times X ,dt\times d\nu)$.\\
Finally we denote by $\bbP_N$ the measure $\bbP_N=\bbP (dw)N_w
(dt,du)$. {One has to remember that  the image of $\bbP\times\nu\times dt$ by
$\crea$ is nothing but $\bbP_N$ whose image by $\anni$ is $\bbP\times\nu\times dt$  (see Lemma 13 in \cite{bouleau-denis})}.
\begin{Th}
The Dirichlet form $(\bbD ,\cE)$ admits a gradient operator that we
denote by $\sharp$ and given by the following formula:
\begin{equation}\label{formulegradient}
\forall F\in\bbD,\quad\ F^\sharp = \int_0^{+\infty}\int_{X\times R}
\anni((\crea F )^\flat)\, dN\odot \rho\in
L^2(\mathbb{P}\times\hat{\mathbb{P}}) .\end{equation} Formula
{\rm(\ref{formulegradient})} is justified by the following decomposition:
$$F\in\mathbb{D}\quad\stackrel{\crea-I}{\longmapsto}\quad \varepsilon^+F-F\in\underline{\mathbb{D}}
\quad\stackrel{\anni((.)^\flat)}{\longmapsto}\quad\anni((\varepsilon^+F)^\flat)\in
L^2_0(\mathbb{P}_N\times\rho)\quad\stackrel{d(N\odot\rho)}{\longmapsto}\quad F^\sharp\in
L^2(\mathbb{P}\times\hat{\mathbb{P}})$$
where each operator is continuous on the range of the preceding one
and where $L^2_0 (\bbP_N \times \rho )$ is the closed set of
elements $G$ in $L^2
(\bbP_N \times \rho )$ such that $\int_R G d\rho=0$ $\bbP_N$-a.e.\\
Moreover, we have for all $F\in\bbD$
\begin{equation}\label{formuleocc}\Gamma [F]=\hE(F^\sharp )^2
=\int_0^{+\infty}\int_X \anni(\gamma [\crea F ])\, dN,\end{equation}
where $\hE$ denotes the expectation with respect to probability
$\hbbP$.
\end{Th}
\begin{proof} This is a slight modification of Theorem 17 in
\cite{bouleau-denis}.
\end{proof}
Let us recall  without proof  some
properties of this structure which are quite natural.\\
\begin{Pro}

 If $h\in L^2 (\R^+,dt)\otimes \bbd$, then $\tN
(h)=\int_0^{+\infty}\int_X h(t,u)\tN (ds,du)$ belongs to $\bbD$ and

\item \begin{equation}
\Gamma [\tN (h)]=\int_0^{+\infty}\int_X \gamma
[h(t,\cdot)](u)N(dt,du).\end{equation}
\item \begin{equation}
\left( \tN (h)\right)^{\sharp}=\int_0^{+\infty}\int_{X\times R}
h^\flat (t,u,r)N\odot\rho (dt,du,dr).
\end{equation}
\end{Pro}
\noindent{\it Remark} 1. Theorem 1 gives a method for obtaining $\Gamma[F]$ for $F\in \mathbb{D}$ or $F\in\mathbb{D}^n $, then with the hypotheses giving (EID) it suffices to prove $\mbox{det }\Gamma[F]>0$ $\mathbb{P}$-a.s.  to assert that $F$ has a density on $\mathbb{R}^n$.
 Let us mention a stronger condition which may be also usefull in some applications. By the following lemma that we leave to the reader
\begin{Le}
Let $M_\alpha$ be random symmetric positive matrices and $\mu(d\alpha)$ a random positive measure. Then
$\{\mbox{\rm det}\int M_\alpha \mu(d\alpha)=0\}\subset\{\int{\mbox{\rm det}}M_\alpha\mu(d\alpha)=0\},$
\end{Le}
\noindent it is enough to have $\int\mbox{det }\varepsilon^-(\gamma[\varepsilon^+F])dN>0$ $\mathbb{P}$-a.s. hence enough that $\mbox{det }\varepsilon^-(\gamma[\varepsilon^+F])$ be $>0$ $\mathbb{P}_N$-a.e. We obtain finally, by lemma 13 of \cite{bouleau-denis}, that a sufficient condition for the density of $F$ is $\mbox{det }\gamma[\varepsilon^+F]>0\;$ $\mathbb{P}\times\nu\times dt$-a.e. (or  equivalently that the components of the vector $(\varepsilon^+F)^\flat$ be $\mathbb{P}\times\nu\times dt$-a.e. linearly independent in $L^2(\rho)$ ).
\subsection{The SDE we consider.}
Let $d\in\Ne$, we consider the following SDE :
\begin{equation}\label{eq}
X_t =x+\int_0^t \int_X c(s,X_{s^-},u)\tN (ds,du)+\int_0^t \sigma
(s,X_{s^-})dZ_s
\end{equation}
where $x\in\R^d$,  $c:\R^+
\times \R^d \times X\rightarrow\R^d$ and $\sigma:\R^+ \times \R^d \rightarrow\R^{d\times n}$ satisfy the set of hypotheses below denoted (R).\\

\noindent{\it Comment.} With respect to the classical form of the SDE's related to Markov processes with jumps, as e.g. Ikeda-Watanabe \cite{ikeda-watanabe} Chap IV \S 9, we put the Brownian part  in the semi-martingale $Z$. Let us emphasize that the Malliavin calculus that we construct, does not concern $Z$ but only $N$. In fact $Z$ could be replaced by a more general random measure with hypotheses assuring existence and uniqueness of the solution. In most applications we have in mind there is no Brownian motion at all, since otherwise this induces strong regularity properties and classical Malliavin calculus applies.\\

\noindent\underline{Assumption (R)}: For simplicity, we fix all along this
article a finite terminal time $T>0$. We assume that $N$ and $Z$ are as explained in \S 2.1. We suppose equation (\ref{eq}) is such that :\\

\noindent 1. There exists $\eta \in L^2 (X,\nu)$  such that:

a) for all $t\in [0,T]$ and $u\in X$, $c(t,\cdot ,u)$ is
differentiable with continuous derivative and
\[\forall u\in X,\  \sup_{t\in [0,T], x\in\R^d} | D_x c(t,x
,u)|\leq \eta (u),\]
\indent b) $\forall (t,u)\in [0,T]\times U,\ |c(t,0,u)|\leq \eta (u)$,

c)  for all $t\in [0,T]$ and $x\in \R^d$, $c(t,x,\cdot )\in\bbd$
and
\[  \sup_{t\in [0,T], x\in\R^d} \gamma [c(t,x,\cdot)](u)\leq \eta(u)
,\]
\indent d) for all $t\in [0,T] $, all $x\in \R^d$ and $u\in X$,\ the matrix $I+D_x
c(t,x,u)$ is invertible and
\[ \sup_{t\in [0,T], x\in \R^d} \left|\left(I+D_x c(t,x,u)\right)^{-1}\right|\leq \eta
(u).\]
2.  For all $t\in [0,T]$ , $\sigma(t,\cdot )$ is
differentiable with continuous derivative and
\[  \sup_{t\in [0,T], x\in\R^d} | D_x \sigma(t,x)|<+\infty.\]
3. As a consequence of  hypotheses 1. \!and 2. \!above, it is well known that equation \eqref{eq} admits a unique solution $X$ such that
$ \E [\sup_{t\in [0,T]} |X_t |^2 ]<+\infty$.
 We suppose that for all $t\in [0,T]$, the matrix $(I+\sum_{j=1}^n D_x \sigma_{\cdot ,j} (t, X_{t^-})\Delta Z_t^j )$ is invertible and its inverse is bounded by a deterministic constant uniformly with respect to $t\in [0,T]$.\\

\noindent{\it Remark} 2. Assumption (R.3) is satisfied if for example one assumes that there exists a constant $a>0$ such that for all $t\in [0,t]$, all $x\in X$ and all $j\in \{ 1,\cdots , n\}$
\[ | \Delta Z_t^j |\leq a \makebox{ and } |D_x\sigma_{\cdot ,j}(t,x)|\leq\frac{1}{na}.\] But it may be verified also without these inequalities. For example when $Z$ is a L\'evy process, because of the independence of the jumps to the past, the invertibility of $(I+\sum_{j=1}^n D_x \sigma_{\cdot ,j} (t, X_{t^-})\Delta Z_t^j )$ is guaranteed in the case $d=n$ by the assumption that the L\'evy measure has a density on $\mathbb{R}^d$ and it remains only to verify the bound of the inverse.
\subsection{Spaces of processes }
 We denote by $\cP$ the predictable sigma-field on
$[0,T]\times \Omega$ and we define the following sets of processes:
\begin{itemize}
\item $\cH$ : the set of real valued  processes $(X_t )_{t\in [0,T]}$, defined on
$(\Omega ,\cA,\bbP)$, which belong to $L^2 ([0,T]\times \Omega)$.
\item $\LP$ : the set of predictable processes in $\cH$.
\item $\LD$ : the set of real valued  processes $(H_t )_{t\in
[0,T]}$, which belong to $L^2 ([0,T]; \bbD)$ i.e. such that
\[ \| H\|^2_{\LD} =\E [ \int_0^T |H_t |^2 dt]+\int_0^T \cE (H_t )dt
<+\infty.\]
\item $\LDP$ : the subvector space of predictable processes in $\LD$.
\item $\LDdP$ : the set of real valued processes $H$ defined on $[0,T]\times
\Omega \times X$ which are predictable and belong to $L^2
([0,T];\bbD\otimes \bbd)$ i.e. such that
\[ \| H\|^2_{\LDdP} =\E [ \int_0^T \int_X |H_t |^2 \nu(du)dt]
+\int_0^T\int_X \cE (H_t(\cdot ,u)) \nu(du)dt+\E [ \int_0^T
e(H_t)dt] <+\infty.\]
\end{itemize}
We define $\LDP^0$ to be the set of elementary processes
in $\LDP$ of the form
\[ G_t
(w)=\sum_{i=0}^{m-1} F_i (w)\1_{]t_i ,t_{i+1}]}(t),\] where
$m\in\Ne$, $0\leq t_0 \leq\cdots t_m\leq T$ and for all $i$, $F_i \in\bbD$ and is $\cA_{t_i}$-measurable.\\
We also consider $\LDdP^0$, the set of elementary processes in
$\LDdP$: $H$ belongs to $\LDdP^0$ if and only if \[ H_t
(w,u)=\sum_{i=0}^{m-1} F_i (w)\1_{]t_i ,t_{i+1}]}(t)g_i (u),\] where for all $i\in \{
0,\cdots ,m-1\}$, $F_i \in\bbD$ and is $\cA_{t_i}$-measurable and
$g_i \in \bbd$.\\
{ The proof of the following lemma is straightforward:}
\begin{Le}  $\LDP^0$ is dense in $\LDP$ and $\LDdP^0$ is dense in $\LDdP$.

 Let $Y\in\LDP$, $c$ and $\sigma$ satisfying
the set of assumptions {\rm(R)}, then
\[ (t,w,u)\longrightarrow c(t,Y_t (w),u)\in\LDdP\]
\[ (t,w)\longrightarrow \sigma (t,Y_t (w))\in \LDP .\]
\end{Le}
\noindent\underline{Notation}: We shall consider $\R^d$-valued processes, so
 $\cH^d$,$\LD^d$,$\ldots$, will denote the spaces of
$\R^d$-valued processes such that each coordinate belongs
respectively to $\cH$,$\LD$,$\ldots$, equipped with the
standard  norm of the product topology.\\
The above spaces  yield the following results allowing to perform the \'Emile Picard iteration procedure with respect to the Dirichlet norm :

\begin{Pro}{\label{Estimees} } Let $H\in \LDdP$ and $G\in \LDP^n$, then:
\begin{enumerate}
\item The process
$$ \forall t\in [0,T], \ X_t =\int_0^t\int_X H(s,w,u)\tN (ds,du)$$
is a square integrable martingale which belongs to $\LD$ and such
that the process $X^- =(X_{t^-})_{t\in [0,T]}$ belongs to $\LDP$.
The gradient operator satisfies for all $t\in [0,T]$:
\begin{equation}\label{Derive1}X_t^\sharp (w,\hw)=\int_0^t\int_X H^\sharp
(s,w,u,\hw ) d\tN (ds,du)+\\\int_0^t\int_{X\times R} H^{\flat}
(s,w,u,r)N\odot\rho (ds,du,dr).\end{equation}
Moreover
\begin{equation}\label{l1}
\forall t\in [0,T],\ \ \parallel X_t \parallel_{\bbD}\leq\sqrt2
\parallel H\parallel_{\LDdP}\qquad\qquad\quad
\end{equation}
\begin{equation}\label{l2}
\quad\parallel X\parallel_{\LD}\leq\sqrt{2T}\parallel
H\parallel_{\LDdP}
\end{equation}
and
\begin{equation}\label{l3}
\parallel X^-\parallel_{\LDP}\leq\sqrt{2T}\parallel
H\parallel_{\LDdP}.
\end{equation}
\item The process
\[ \forall t\in [0,T], \ Y_t =\int_0^t G(s,w)dZ_s\]
is a square integrable semimartingale which belongs to $\LD$,
$Y^-=(Y_{t^-})_{t\in [0,T]}$ belongs to $\LDP$ and
\begin{equation}\label{Derive2}\forall t\in [0,T],\ Y_t^\sharp (w,\hw)=\int_0^t G^\sharp
(s,w,\hw ) dZ_s.\end{equation} We also have the following estimates:
\[\forall t\in [0,T],\ \parallel Y_t \parallel_{\bbD}\leq C \parallel
G\parallel_{\LDP^n} \makebox{ and } \parallel Y_{t^-}
\parallel_{\bbD}\leq C\parallel G\parallel_{\LDP^n},\]
\[ \parallel Y\parallel_{\LD}\leq C\sqrt{T}\parallel
G\parallel_{\LDP^n}\makebox{ and } \parallel Y^-
\parallel_{\LDP}\leq C\sqrt{T}\parallel G\parallel_{\LDP^n}.\]
\end{enumerate}
\end{Pro}
\begin{proof}
1. Assume first that $H\in\LDdP^0$,
\[ H
(w,s,u)=\sum_{i=0}^{m-1} F_i (w)\1_{]t_i ,t_{i+1}]}(s)g_i (u).\]
Then, $X$ is a square integrable martingale and we have: \[ \forall t\in [0,T], X_t
=\sum_i F_i \tN (\1_{]t_i \wedge t,t_{i+1}\wedge t]}\cdot g_i ).\] So, $X_t$ belongs
to $\bbD$ and thanks to the functional calculus :
\begin{eqnarray*}
X_t^\sharp &=& \sum_i \left( F_i^\sharp\cdot\tN (\1_{]t_i \wedge t
,t_{i+1}\wedge t]}\cdot g_i )+F_i \cdot\left(\tN (\1_{]t_i \wedge t,t_{i+1}\wedge t]}\cdot
g_i )\right)^\sharp \right)\\
&=&\sum_i \left( F_i^\sharp\cdot\tN (\1_{]t_i \wedge t,t_{i+1}\wedge t]}\cdot g_i
)+F_i \cdot\int_{]t_i \wedge t ,t_{i+1}\wedge t]}\int_{X\times R} g_i^\flat
(u,r)N\odot\rho (ds,du,dr)\right)\\
&=&\int_0^t\!\int_X H^\sharp (w,s,u,\hw ) d\tN
(ds,du)+\int_0^t\!\int_{X\times R} H^{\flat} (w,s,u,r)N\odot\rho
(ds,du,dr).
\end{eqnarray*}
This yields:
\begin{eqnarray*}
\cE(X_t)&=&\frac{1}{2}\E\hE|X^\sharp |^2\\
&\leq&  \E\left[\int_0^t\!\int_X \hE\left[|H^\sharp (s,u,\hw
)|^2\right]\nu(du)ds\right]+\E\left[\int_0^t\!\int_{X\times R}
|H^{\flat}(s,u,r)|^2
\nu(du)\rho(dr)ds\right]\\
&=&2\left(\int_0^t\int_X \cE (H(s ,u) \nu(du))ds+\E [ \int_0^t e(H
(s,\cdot))ds]\right),
\end{eqnarray*}
so we obtain inequalities (\ref{l1}) and (\ref{l2}) in this case and
then for general $H\in\LDdP$ by
density. \\
As both $X$ and $X^\sharp$ are stochastic integrals w.r.t. $\!$random
measure, we deduce that $\lim_{s\rightarrow t, s<t} X_t$ exists in
$L^2 (\Omega, \bbP)$, $\lim_{s\rightarrow t, s<t} X^\sharp_t$ exists
in $L^2 (\Omega\times \hat{\Om},\bbP\times\hbbP)$ so that, by
standard arguments, $\lim_{s\rightarrow t, s<t} X_t$ exists in
$\bbD$ and as we have
\[\forall s<t,\ \parallel X_s \parallel_{\bbD}\leq \sqrt2\parallel
H\parallel_{\LDdP},\] what yields this part of the proposition.

\noindent 2.  For the second part, consider first $G\in (\LDP^0)^n $
\[ G_t
(w)=\sum_{i=0}^{m-1} F_i (w)\1_{]t_i ,t_{i+1}]}(t).\]
We have
\begin{eqnarray*}
Y_t& =&\sum_{i=0}^{m-1} F_i  (Z_{t_{i+1}\wedge t}-Z_{t_i \wedge t}),\\
Y_t^\sharp& =&\sum_{i=0}^{m-1} F_i^\sharp  (Z_{t_{i+1}\wedge t}-Z_{t_i \wedge t}).
\end{eqnarray*}
Thanks to the bounds \eqref{CondCrochet}, we easily obtain
\begin{eqnarray*}
\E [|Y_t |^2 ] & \leq& C\mathbb{E}\int_0^t |G_s |^2\,
ds,
\end{eqnarray*}
and
$$
\cE (Y_t )=\frac{1}{2}\E\hE [|Y_t^\sharp |^2]
\leq\frac{C^2}{2} \hE\E \int_0^t \,|G^\sharp_s |^2 ds
=C^2 \int_0^t \sum_{j=1}^n \cE (G^j_s )ds.$$

It is now easy to conclude, using a density argument, similarly to
the proof of the first part.
\end{proof}
\subsection{Preliminary results on the solution of equation (\ref{eq}).}
\begin{Pro} The equation {\rm(\ref{eq})} admits a unique solution $X$ in
$\LD^d$. Moreover, the gradient of $X$ satisfies:
\begin{eqnarray*}
X_t^\sharp &=&\int_0^t\int_U D_x c(s,X_{s-},u)\cdot
X^\sharp_{s-}\tN (ds,du)+\int_0^t\int_{X\times R}
c^{\flat}(s,X_{s-},u,r)N\odot\rho (ds,du,dr)\\&&+\int_0^t D_x
\sigma (s,X_{s-})\cdot X^\sharp_{s-} dZ_s
\end{eqnarray*}
\end{Pro}
\begin{proof} { We define inductively a sequence $(X^r)$ of
$\R^d$-valued semimartingales by $X^0 =x$ and
\[\forall r\in\N,\ \forall t\in [0,T],\ X^{r+1}_t =x+
\int_0^t \int_X c(s,X^r_{s^-},u)\tN (ds,du)+\int_0^t \sigma
(s,X^r_{s^-})dZ_s.\] As a consequence of Proposition
\ref{Estimees}, it is clear that for all $r$, $X^r$ belongs to
$\LD^d$ and that we have $\forall t\in [0,T]$
 \begin{eqnarray*}
X_t^{r+1 ,\sharp}& =&\int_0^t\int_U D_x c(s,X^r_{s-},u)\cdot
X^{r,\sharp}_{s-}\tN (ds,du)+\int_0^t\int_{U\times R}
c^{\flat}(s,X^r_{s-},u,r)N\odot\rho (ds,du,dr)\\&&+\int_0^t D_x
\sigma (s,X^r_{s-})\cdot X^{r,\sharp}_{s-} dZ_s.\end{eqnarray*}
This  is the iteration procedure due to \'Emile Picard and it
is well-known that \begin{equation}\label{Picard} \lim_{r\rightarrow +\infty}E[\sup_{t\in
[0,T]} |X_t -X^r_t |^2] =0.\end{equation} Moreover, thanks to the hypotheses
we made on the coefficients, it is easily seen that there exists a
constant $\kappa$ such that for all $r\in\Ne$ and all $t\in [0,T]$
\begin{eqnarray*}
\E\hE\left[ |X_t^{r+1 ,\sharp}|^2\right] &\leq&\kappa \left( 1+\int_0^t \E\hE
\left[|X^{r,\sharp}_{s-}|^2 \right]ds \right)\end{eqnarray*} so that by
induction we deduce
\[ \forall r\in\N ,\ \forall t\in [0,T] ,\ \E\hE \left[|X_t^{r ,\sharp}|^2\right]\leq \kappa e^{\kappa t}.\]
Hence, the sequence $(X^r )$ is bounded in $\LD^d$ which is an Hilbert space. Therefore, there is a sequence of convex combinations of
$X^r$ which converges to a process $Y\in \LD^d$. But, by \eqref{Picard} we a priori know that $X^r$ tends to $X$ in $L^2 ([0,T];\R^d )$ so that $Y$ is nothing but  $X$. This proves that $X$ belongs to $\LD^d$ and the relation satisfied by the gradient is consequence of relations \eqref{Derive1} and \eqref{Derive2}.
}\end{proof}
 We can now explicit a formula for the carré du champ
operator of $X_t$, using the linear equation satisfied by
$X^{\sharp}$. The obtained formula (Theorem \ref{OCC} below) is
known for long time (cf \cite{bichteler-gravereaux-jacod}) but
established here under much weaker regularity assumptions similar
to those of \cite{coquio}. Let us emphasize that we obtain this
formula without the intensity of the Poisson measure being the
Lebesgue measure as supposed by these authors.

Let us define the  $\R^{d\times d}$-valued process $U_s$ by
$$dU_s=\sum_{j=1}^nD_x\sigma_{.,j}(s,X_{s-})dZ_s^j.$$
Then the following $\R^{d\times d}$-valued
process is the derivative of the flow generated
by $X$:
\begin{eqnarray*}
K_t &=& I+\int_0^t\int_X D_x c(s,X_{s-} ,u)K_{s-} \tN (ds ,du)
+\int_0^t dU_sK_{s-}
\end{eqnarray*}

 \noindent Under our
hypotheses, for all $t\geq 0$, the matrix $K_t$ is invertible as a consequence of the following proposition which extends classical formulas about linear equations (e.g. \cite{protter} Chap V \S 9 Thm 52).

\begin{Pro}\label{calcul-matriciel} Let $\Sigma_t$ be a $d\times d$-matrix semimartingale such that $I+\Delta\Sigma_t$  is invertible $\forall t$ a.s.
Let $K_t$ be the solution of
\begin{equation}\label{K}K_t=I+\int_0^td\Sigma_sK_{s-}.\end{equation}
a) Then $K_t$ is invertible and its inverse $\bK _t$ satisfies
\begin{equation}\label{K-1}\bK_t=I-\int_0^t\bK _{s-}d\Sigma_s+\sum_{s\leq t}\bK _{s-}(\Delta\Sigma_s)^2(I+\Delta\Sigma_s)^{-1}+\int_0^t\bK _sd<\Sigma^c,\Sigma^c>_s.\end{equation}
b) The $d\times1$-solution of
$V_t=V_\alpha+\int_{]\alpha}^td\Sigma_sV_{s-}$
is given by $V_t=K_t\bK_\alpha V_\alpha.$

\noindent c) Let $R_t$ be a $d\times1$-semimartingale. The solution of
$S_t=R_t+\int_0^td\Sigma_sS_{s-}$ is given by
$$S_t=R_0+K_t[\int_0^t\bK _{s-}dR_s-\sum_{s\leq t}\bK _{s-}\Delta\Sigma_s(I+\Delta\Sigma_s)^{-1}\Delta R_s-\int_0^t\bK _sd<\Sigma^c,\Sigma^c>_s].$$

\end{Pro}
\begin{proof} By Ito's formula the solutions to (\ref{K}) and (\ref{K-1}) satisfy
$$\begin{array}{rl}dK_t\bK _t=&d\Sigma_tK_{t-}\bK _{t-}-K_{t-}\bK _{t-}d\Sigma_t\\
&+[K_{t-}\bK _{t-}(\Delta\Sigma_t)^2-\Delta\Sigma_tK_{t-}\bK _{t-}\Delta\Sigma_t](I+\Delta\Sigma_t)^{-1}\\
&+K_{t-}\bK _{t-}d<\Sigma^c,\Sigma^c>_t-<d\Sigma^cK_-,\bK Ñ_-d\Sigma^c>_t
\end{array}$$
this equation is Lipschitz with respect to the unknown $K\bK $ and is verified by the identity matrix, what gives the first assertion of the proposition. The remaining is similar.
\end{proof}
\noindent Here $\bK_t =(K_t )^{-1}$ satisfies:
\begin{eqnarray*}
\bK_t &=& I -\int_0^t\int_X \bK_{s-} (I+D_x c(s,X_{s-} ,u))^{-1} {
D_x c(s,X_{s-} ,u)}\tN (ds ,du)
\\&&-\int_0^t \bK_{s-} dU_s+\sum_{s\leq t}\bK_{s-}(\Delta U_s)^2(I+\Delta U_s)^{-1}+\int_0^t\bK_s d<U^c,U^c>_s.
\end{eqnarray*}
\subsection{Obtaining the Malliavin matrix thanks to the lent particle method.}
\begin{Th}{\label{OCC}}
For all $t\in [0,T]$,
\begin{eqnarray*}
\Gamma [X_t ]&=&K_t  \int_0^t \int_X\bK_{s} \gamma[c(s
,X_{s-} ,\cdot )]\bK_{s}^{\ast}\, N (ds, du)
K_t^{\ast},
\end{eqnarray*}
 where for any matrix $D$, $D^{\ast}$ denotes its transpose.
\end{Th}
\begin{proof} We apply the lent particle method:

Let $(\alpha ,u)\in [0,T]\times X$. We put
$ X_t\alu =\crea_{(\alpha ,u)} X_t.$
We have,
\begin{eqnarray*}
X_t\alu &=&x+\int_0^{\alpha} \int_X c(s,X_{s^-}\alu,u')\tN
(ds,du')+\int_0^{\alpha} \sigma(s,X_{s^-}\alu) dZ_s+c(\alpha
,X_{\alpha^-}\alu ,u)\\&&+\int_{]\alpha, t]} \int_X
c(s,X_{s^-}\alu,u')\tN (ds,du')+\int_{]\alpha, t]} \sigma
(s,X_{s^-}\alu)dZ_s.
\end{eqnarray*}
Let us remark that $X_t\alu =X_t $ if $t<\alpha$ so that, taking
the gradient with respect to the variable $u$, we obtain:
\begin{eqnarray*}
(X_t\alu )^\flat&=&(c (\alpha ,X_{\alpha^-}\alu ,u))^\flat
+\int_{]\alpha, t]} \int_X D_x c(s,X_{s^-}\alu,u')\cdot
(\xs\alu)^\flat \tN (ds,du')\\&&+\int_{]\alpha, t]}
D_x\sigma(s,X_{s^-}\alu)\cdot (\xs\alu)^\flat dZ_s.
\end{eqnarray*}
 Let us now introduce
the process $K_t\alu = \crea_{(\alpha ,u)}(K_t)$ which satisfies
the following SDE:
$$
K_t\alu  = I+\int_0^t\int_X D_x c(s,\xs\alu ,u') K_{s-}\alu \tN (ds ,du')
+\int_0^t dU_s^{(\alpha,u)}K_{s-}^{(\alpha,u)}
$$
and its inverse $\bK_t\alu = (K_t\alu)^{-1}$.
Then, using the flow property, (Prop (\ref{calcul-matriciel})b), we have:
\[\forall t\geq 0 ,\ (X_t\alu )^\flat = K_t\alu  \bK_{\alpha}\alu
(c(\alpha ,X_{\alpha^-} ,u))^\flat .\]
(Let us note that at this stage Remark 1 of Section 2.2 could be used to get a sufficient condition of density of $X_t$.)

\noindent Keeping in mind that the measures $\mathbb{P}\times\nu$ and $\mathbb{P}_N$ are mutually singular, let us emphasize that this result is an equality $\mathbb{P}\times\nu(du)\times\rho$-a.e. (cf Thm 1 above or \cite{bouleau-denis} Prop 18).
Now, we calculate the carr\'e
du champ and then we take back the particle:
\[\forall t\geq 0 ,\ \anni_{(\alpha ,u)}\gamma [ (X_t\alu )] = K_t  \bK_{\alpha}
\gamma[c(\alpha ,X_{\alpha^-} ,\cdot )]
\bK_{\alpha}^{\ast} K_t^{\ast}
\] which is an equality $\mathbb{P}_N$-a.e. (cf Thm 1 or \cite{bouleau-denis} Prop 18).

Finally integrating with respect to $N$ we get
\begin{eqnarray*}
\forall t\geq 0 ,\ \Gamma [X_t]& =& K_t  \int_0^t
\int_X\bK_{s} \gamma[c(s ,X_{s^-} ,\cdot )](u)
\bK_{s}^{\ast}N (ds, du) K_t^{\ast}\end{eqnarray*} which
ends the proof.
\end{proof}
\noindent{\it Remark} 3. We could also write $\Gamma [X_t]$ as
\[K_t  \int_0^t\!\! \int_X\bK_{s^-} (I+D_x
c(s,\xs ,u))^{-1}  \gamma[c(s ,X_{s^-} ,\cdot )](u)
(I+(D_x c(s,\xs ,u))^\ast)^{-1} \bK_{s^-}^{\ast}\, N (ds,
du) K_t^{\ast}.\]
\noindent{\it Remark} 4. With the assumptions (R) that we have taken on the SDE  (\ref{eq}),  $X_t$ is in $\mathbb{D}$ and $\Gamma [X_t]$ given by the above formulae  is in $L^1(\mathbb{P})$. Now, let us recall that there exists a powerfull Borelian localization procedure in any local Dirichlet structure with carr\'e du champ (cf \cite{bouleau-hirsch2} Chap. I \S7.1). In practice, assumptions (R) may be lightened in such a way that significant positive quantities be only finite almost everywhere, then the above formulae are still true with  $X_t\in\mathbb{D}_{loc}$ and  $\Gamma [X_t]$ finite a.s. and property (EID) applies as well.

\section{Some applications.}
{
\subsection{The regular case}
 This is the case where we assume that  $X$ is a topological space and that coefficient $c$ is regular with respect to the jumps size. More precisely, we have the following

 \begin{Pro} Assume that X is a topological space, that the intensity measure $ds\times\nu$ of $N$ is such that $\nu$ has an infinite mass near some point $u_0$ in $X$. If the matrix $(s,y,u)\rightarrow\gamma[c(s,y,\cdot)](u)$ is continuous on a neighborhood of $(0,x,u_0)$ and invertible at $(0,x,u_0)$, then the solution $X_t$ of {\rm(\ref{eq})} has a density for all $t\in ]0,T]$.
 \end{Pro}
 \begin{proof} Let us fix $t\in ]0,T]$. As $\nu$ has infinite mass near $u_0$, as $X$ is right continuous and $\gamma [c]$ continuous, $N$ admits almost surely  a jump at time $s\in ]0, t]$ with size $u\in X$  such that $\gamma [c(s, X_{s^-},\cdot)](u)$ is invertible. As a consequence,
 \[ \Gamma [X_t ]\geq K_t  \bK_{s} \gamma[c(s ,X_{s^-} ,\cdot )](u)
\bK_{s}^{\ast} K_t^{\ast},\] for the relation order  in the set of
non-negative symmetric matrixes. As $\Gamma [X_t ]$ dominates an
invertible matrix, it is also invertible and this permits to
conclude.
\end{proof}
}
Now, the method can yield existence of density for $X_t$ even when the matrix $\gamma[c(s,x,\cdot)]$ is everywhere singular. A example of such a situation was given in \cite{bouleau-denis} \S 5.3.\\
 {We now turn out to study other degenerated examples.}

\subsection{L\'evy's stochastic area.}
This example will show that the method can detect densities even when both the matrix $\gamma[c]$ is non invertible and the  L\'evy measure singular.

Let $X(t)=(X_1(t), X_2(t))$ be a  L\'evy process with values in $\mathbb{R}^2$ with  L\'evy measure $\sigma$. We suppose  that the hypotheses of the method are fulfilled, we shall explicit this later on.

Let us consider for the moment a general gradient on the bottom space :
$$f^\flat=f_1^\prime\xi_1+f_2^\prime\xi_2$$
{ where $f'_i =\frac{\partial f}{\partial x_i}, $ and
$\xi_1$, $\xi_2$ are functions defined on $\R^2 \times R$ which
satisfy: $\int_R \xi_1 (\cdot , r)\rho (dr)=\int_R \xi_2 (\cdot ,r
)\rho (dr)=0$, $\int_R \xi_1^2 (x_1 ,x_2 ,r)\rho
(dr)=\alpha_{11}(x_1,x_2)$, $\int_R \xi_1 (x_1 ,x_2 ,r)\xi_2 (x_1
,x_2 ,r)\rho (dr)=\alpha_{12}(x_1,x_2)$, $\int_R \xi_2^2 (x_1 ,x_2
)\rho (dr)=\alpha_{22}(x_1,x_2)$}, so that
$$\gamma[f]=\alpha_{11}f^{\prime 2}_1+2\alpha_{12}f^\prime_1f^\prime_2+\alpha_{22}f_2^{\prime2}.$$
Let be
$$V=(X_1(t),X_2(t),\int_0^tX_1(s_-)dX_2(s)-\int_0^tX_2(s_-)dX_1(s)).$$
We have { for $0<\alpha < t$ and $x=(x_1 ,x_2
)\in\R^2,$}
$$\varepsilon^+_{(\alpha ,x)}V=V+(x_1,x_2, X_1(\alpha_-)x_2+x_1(X_2(t)-X_2(\alpha))-X_2(\alpha_-)x_1-x_2(X_1(t)-X_1(\alpha))$$
$$\qquad=V+(x_1,x_2,x_1(X_2(t)-2X_2(\alpha))-x_2(X_1(t)-2X_1(\alpha)))$$
because $\varepsilon^+V$ is defined
$\mathbb{P}\times\nu{ \times d\alpha}$-a.e. and $\nu\times d\alpha$ is
diffuse, so
$$(\varepsilon^+V)^\flat=(\xi_1,\xi_2,\xi_1(X_2(t)-2X_2(\alpha))-\xi_2(X_1(t)-2X_1(\alpha)))$$
and
$$
\gamma[\varepsilon^+V]=\left(
\begin{array}{ccc}
\alpha_{11}&\alpha_{12}&A\alpha_{11}-B\alpha_{12}\\
\alpha_{12}&\alpha_{22}&A\alpha_{12}-B\alpha_{22}\\
A\alpha_{11}-B\alpha_{12}&A\alpha_{12}-B\alpha_{22}&A^2\alpha_{11}-2AB\alpha_{12}+B^2\alpha_{22}
\end{array}
\right)
$$
denoting
$A=(X_2(t)-2X_2(\alpha))$ et $B=(X_1(t)-2X_1(\alpha))$.

The operator $\varepsilon^-$ gives a functional defined
$\mathbb{P}_N$-a.e. \! so that for example $$\varepsilon^-_{(\alpha,x_1
,x_2)}(X(t))=X(t)-\Delta X_{\alpha}\qquad\mathbb{P}_N(d\alpha dx_1dx_2)\mbox{-a.e.}$$ This yields
$$\varepsilon^-A=X_2(t)-\Delta X_2(\alpha)-2X_2(\alpha_-)\qquad\mbox{let us denote it }\tilde{A}$$
$$\varepsilon^-B=X_1(t)-\Delta X_1(\alpha)-2X_1(\alpha_-)\qquad\mbox{let us denote it }\tilde{B}$$
and eventually
$$
\Gamma[V]=\sum_{\alpha\leq t}\left(
\begin{array}{ccc}
\alpha_{11}(\Delta X_\alpha)&\alpha_{12}(\Delta X_\alpha)&\tilde{A}\alpha_{11}(\Delta X_\alpha)-\tilde{B}\alpha_{12}(\Delta X_\alpha)\\
\sim&\alpha_{22}(\Delta X_\alpha)&\tilde{A}\alpha_{12}(\Delta X_\alpha)-\tilde{B}\alpha_{22}(\Delta X_\alpha)\\
\sim&\sim&\tilde{A}^2\alpha_{11}(\Delta X_\alpha)-2\tilde{A}\tilde{B}\alpha_{12}(\Delta X_\alpha)+\tilde{B}^2\alpha_{22}(\Delta X_\alpha)
\end{array}
\right)
$$ the symbol $\sim$ denoting the symmetry of the matrix.

\subsubsection{3.2.1. First case.} Let us consider the case $\alpha_{12}=0$. We are in this case if $\nu$ possesses a density satisfying our hypotheses which are fulfilled as soon as we assume those of Lemma \ref{HypBC}, and under these assumptions, with same notation,  we may choose $\alpha_{11}=\alpha_{22}=((x_1^2+x_2^2)\frac{\psi (x)}{k(x)})$. We have

$$\Gamma[V]=\sum_{\alpha\leq t}|\Delta X_\alpha|^2 \frac{\psi (\Delta X_\alpha)}{k(\Delta X_\alpha)}\left(
\begin{array}{ccc}
1&0&\tilde{A}\\
0&0&0\\
\tilde{A}&0&\tilde{A}^2
\end{array}
\right)
+|\Delta X_\alpha|^2 \frac{\psi (\Delta X_\alpha)}{k(\Delta X_\alpha)}\left(
\begin{array}{ccc}
0&0&0\\
0&1&\tilde{B}\\
0&\tilde{B}&\tilde{B}^2
\end{array}
\right)
$$

\noindent Hence   $X$ has a density if the dimension of the vector space spanned by
$$\left(\left(
\begin{array}{c}
1\\
0\\
X_2(t)-\Delta X_2(\alpha)-2X_2(\alpha-)
\end{array}
\right),
\left(
\begin{array}{c}
0\\
1\\
X_1(t)-\Delta X_1(\alpha)-2X_1(\alpha-)
\end{array}
\right),
\alpha\in JT\right)$$
is equal to 3, where $JT=\{ \alpha \in [0,t],\ \Delta X_{\alpha}\in O\}$.

Let us suppose that $\nu (O)=+\infty$. Then one of the projections of $O$ on the axes has an infinite mass. Let us suppose it is that of  $X_1$.

The process $X_1(t)-\Delta X_1(\alpha)-2X_1(\alpha-)=X_1(t)-X_1(\alpha)-X_1(\alpha-)$ cannot remain constant for $\alpha\in JT$ since the L\'evy measure of $X_1$ is infinite, and therefore  $JT$ is dense in $\mathbb{R}_+$ and each point of $JT$ is the limit of an increasing sequence  $(\alpha_k)$ of points of $JT$ such that $|\Delta X_1(\alpha_k)|\rightarrow 0$. The constancy of $X_1(t)-\Delta X_1(\alpha)-2X_1(\alpha-)$ would imply that of $-2X_1(\alpha-)$ and there will be no jumps. Thus,   $V$ has a density  if {\it the L\'evy measure, $\nu$, of $X$ satisfies
hypotheses of Lemma \ref{HypBC}  } and $\nu (O)=+\infty$.\\

\noindent{\it Example 1.} Let us take the  L\'evy measure of $(X_1,X_2)$ expressed in polar coordinates as
$$\nu(d\rho,d\theta)=g(\theta)d\theta.1_{]0,1[}(\rho)\frac{d\rho}{\rho}$$
with $g$ locally bounded and such that it dominates a continuous and positive function near $0$.
Then $V=(X_1(t),X_2(t), \int_0^tX_1(s_-)dX_2(s)-\int_0^tX_2(s_-)dX_1(s))$ has a density (and condition (0.4) of \cite{cancelier-chemin} or of \cite{picard} prop1.1 are not fulfilled).\\

\subsubsection{3.2.2. Second case.} Let us suppose $\xi_2=\lambda(x_1,x_2)\xi_1$.
We are in this case if the measure $\nu$ is carried by a
graph in $\mathbb{R}^2$ and image of a measure on $\mathbb{R}$. For instance if $X_2$ is taken to be $[X_1]$ the L\'evy measure is carried by the graph $x_2=x_1^2$.

Then
$\alpha_{12}=\lambda\alpha_{11}$ et
$\alpha_{22}=\lambda^2\alpha_{11}$. We have
$$
\Gamma[V]=\sum_{\alpha\leq t}\alpha_{11}(\Delta X_\alpha)\left(
\begin{array}{ccc}
1&\lambda&\tilde{A}-\lambda\tilde{B}\\
\lambda&\lambda^2&\lambda\tilde{A}-\lambda^2\tilde{B}\\
\tilde{A}-\lambda\tilde{B}&\lambda\tilde{A}-\lambda^2\tilde{B}&(\tilde{A}-\lambda\tilde{B})^2
\end{array}
\right).
$$ where $\lambda$ is taken on the jumps : $\lambda(\Delta X_\alpha)$.
If $\alpha_{11}\neq0$, $V$ has a density as soon as
\begin{equation}\label{span}\mbox{dim }\mathcal{L}\left(\left(
\begin{array}{c}
1\\
\lambda\\
\tilde{A}-\lambda\tilde{B}
\end{array}
\right),\quad\alpha\in JT\right)=3\end{equation}

\noindent with $\tilde{A}-\lambda\tilde{B}=-X_2(\alpha_-)+\lambda(\Delta X(\alpha))X_1(\alpha_-)+X_2(t)-X_2(\alpha)-\lambda(\Delta X(\alpha))(X_1(t)-X_1(\alpha)).$\\

In order to study condition (\ref{span}) let us reason on the set
$$A=\{\omega : \mbox{dim }\mathcal{L}\left(\left(
\begin{array}{c}
1\\
\lambda\\
\tilde{A}-\lambda\tilde{B}
\end{array}
\right),\quad\alpha\in JT\right)<3\}
$$
There exist $a,b,c$ (d\'epending on $\omega$) such that $\forall\alpha\in JT$
$$a+b\lambda+c(\tilde{A}-\lambda\tilde{B})=0.$$ If the L\'evy measure of $X$ is infinite,  $JT$ is dense in $\mathbb{R}_+$, each point of $JT$ is limit of an increasing sequence  $(\alpha_k)$ of points of $JT$ such that $|\Delta X(\alpha_k)|\rightarrow 0$. If the function $\lambda$ goes to zero at zero, the process
$a+c(X_2(t)-2X_2(\alpha-))$ vanishes on $JT$.
If the L\'evy measure of $X_2$ is also infinite,  $X_2(\alpha-)$ cannot remain constant when $\alpha$ varies, hence $c=a=b=0$.

Thus,
\noindent{\it choosing $\alpha_{11}=x_1^2\wedge 1$ (or $\alpha_{11}=X_1^2 \frac{\psi (x)}{k(x)}$ if one assumes hypotheses of Lemma \ref{HypBC}), if $\lambda $  tends to zero at zero, if the  L\'evy measure of $X_2$ is  infinite, and such that there exists a bottom structure $(\mathbb{R}^2\backslash\{0\}, \mathcal{B},\sigma, \mathbf{ d},\gamma)$ allowing {\rm(BC)} and {\rm(EID)} on the upper space, then $V$ has a density.}\\

\noindent{\it Example 2.} This applies to $V=(X_1(t),[X_1]_t, \int_0^tX_1(s_-)d[X_1](s)-\int_0^t[X_1](s_-)dX_1(s))$.

 The L\'evy measure of $(X_1,[X_1])$ is carried by the curve  $x_2=x_1^2$. We have $\lambda(x_1,x_2)=2x_1$. We arrive to the sufficient condition :
{\it $V$ has a  density as soon as the L\'evy measure of $X_1$ satisfies hypotheses of Lemma \ref{HypBC} with $\nu (O)=+\infty$}.

\subsection{McKean-Vlasov type equation driven by a L\'evy process.}

The following nonlinear stochastic differential equation
\begin{equation}\label{nonlinear}
\left\{\begin{array}{l}
X_t=X_0+\int_0^t\sigma(X_{s-},P_s)\;dY_s\quad t\in [0,T]\\
\forall s\in[0,T], P_s \mbox{ is the probability law of } X_s
\end{array}\right.
\end{equation}
where $Y$ is a L\'evy process with values in $\mathbb{R}^d$, independent of $X_0$, and $\sigma$ : $\mathbb{R}^k\times\mathcal{P}(\mathbb{R}^k)\mapsto\mathbb{R}^{k\times d}$ where $\mathcal{P}(\mathbb{R}^k)$ denotes the set of probability measures on $\mathbb{R}^k$, generalizes the McKean-Vlasov model. It has been studied by Jourdain, M\'el\'eard and Woyczynski \cite{jourdain} who proved, by a fixed point argument, that equation (\ref{nonlinear}) admits a solution as soon as $\sigma$ is Lipschitz continuous on $\mathbb{R}^k\times\mathcal{P}(\mathbb{R}^k)$ equipped with the product of the canonical metrics on $\mathbb{R}^k$ and a modified Wasserstein metrics on probability measures.

When $Y$ is a one-dimensional L\'evy process and $k=1$, these authors obtained the existence of a density for $X_t$ using a Malliavin calculus in the Bichteler-Gravereaux-Jacod spirit under the assumptions that $\sigma$ does not vanish, admits two bounded derivatives with respect to the first variable, and the L\'evy measure of $Y$ dominates an absolutely continuous measure with $\mathcal{C}^2$-density and infinite mass, and additional technical conditions.

We would like to illustrate the lent particle method by simplifying their proof and lightening the hypotheses.

The clever remark --- evident after a while of reflection --- used by these authors, and usefull for us too, is that as soon existence for equation (\ref{nonlinear}) has been proved, it may be considered for Malliavin calculus as an equation of the form
$$X_t=X_0+\int_0^ta(X_{s-},s)\;dY_s$$
which is a particular case of our present study.

Let us proceed with the following hypotheses :

(i) $a$ is $\mathcal{C}^1\cap Lip$ with respect to the first variable uniformly in $s$ and
$$\sup_{t,x}|(I+D_xa)^{-1}(x,t)|\leq \eta$$

(ii) the L\'evy measure of $Y$ is such that a Dirichlet structure may be chosen such that (H0) and (EID) be fulfilled on the Poisson space (we shall detail this assumption later on).

By the lent particle method we obtain
\begin{equation}\label{Gamma}
\Gamma[X_t]=K_t\left[\sum_{\alpha\in JT}\overline{K}_\alpha a(X_{\alpha-},\alpha)\gamma[j,j^\ast](\Delta Y_\alpha)a^\ast(X_{\alpha-},\alpha)\overline{K}_\alpha^\ast \right]K^\ast_t
\end{equation}
where $JT$ is the random set of jump times of $Y$ before $t$, $\gamma$ is the carr\'e du champ of the bottom structure, $j$ is the identity map on $\mathbb{R}^d$, $K_t$ is the solution of $$K_t=I+\int_0^t\,dZ_s\, K_{s-}$$ where $dZ_s=\sum_{i=1}^dD_xa_{.i}(X_{s-},s)\;dY_s^i$, and $\overline{K}$ the inverse of $K$.

As a consequence of formula (\ref{Gamma}) we can conclude that $X_t$ possesses a density on $\mathbb{R}^k$ under the following hypotheses :

$1^o$) the L\'evy measure of $Y$ satisfies hypotheses of Lemma \ref{HypBC} with $\nu (O)=+\infty$. Then we may choose the operator $\gamma$ to be
$$\gamma[f]=\frac{\psi (x)}{k(x)}\sum_{i=1}^d x_i^2\sum_{i=1}^d(\partial_if)^2\quad \mbox{for }f\in\mathcal{C}^{\infty}_0(\mathbb{R}^d)$$ and the identity map $j$ belongs to $\mathbf{ d}$ and $\gamma[j,j^\ast](x)=\frac{\psi (x)}{k(x)}|x|^2I$. (see \cite{bouleau-denis} for a weaker assumption and the proof of (EID) on the Poisson space).

$2^o$) $a$ satisfies (i), is continuous with respect to the second variable at 0, and such that the matrix $aa^\ast(X_0,0)$ is invertible.\\

\subsection{Stable-like processes.}

The passage between a L\'evy kernel $\nu(t,x,dy)$ with which is expressed the generator of a Markov process with jumps, with standard notation
\begin{equation}\label{L}
\begin{array}{rl}
 Lf(t,x)&=\displaystyle{\frac{1}{2}\sum_{ij=1}^da_{ij}(t,x)\partial^2_{ij}f(x)+\sum_{i=1}^d\partial_if(x)}\\
&\displaystyle{+\int_{\mathbb{R}^d\backslash\{0\}}\left(f(x+y)-f(x)-1_{|y|<1}\sum_{i=1}^d\partial_if(x)y_i\right)}\nu(t,x,dy)
\end{array}
\end{equation}
to the Poisson random measure $N(dt,du)$ to be used for the SDE able to yield  the Markov process as solution
$$dX_t=\sigma(t,X_t)dB_t+b(t,X_t)dt+\int c(t,X_{t-},u)\tilde{N}(dt,du)$$
is theoretically always possible thanks to a result of El Karoui-Lepeltier \cite{elkaroui-lepeltier}. But this general procedure yields for $c(t,x,u)$ a function with few regularity (as the theorem allowing to simulate any probability law on $\mathbb{R}^d$ thanks to a random variable defined on $[0,1]$ equipped with the Lebesgue measure, cf \cite{bouleau-lepingle} Chap. I Thm A.3.1).

The study of the correspondence between $\nu(t,x,dy)$ and the pair $(c(t,x,u), N(dt,du))$ has been deepened by Tsuchiya \cite{tsuchiya} in order to find conditions on $\nu(t,x,dy)$ such  that a function $c(t,x,u)$ may be found verifying the Lipschitz hypotheses guaranteeing the existence of a solution to the SDE, hence of the Markov process with generator (\ref{L}).

He applies this study to the case of so-called {\it stable-like} processes {\it of order $\alpha(x)$} introduced in dimension one by Bass \cite{bass}, whose generator may be symbolically written
\begin{equation}\label{LL}
L=-(-\Delta)^{\frac{\alpha(x)}{2}}
\end{equation} and he obtains the existence of the Markov process for $\alpha$ Lipschitz and such that $0\leq \alpha(x)<2$. (see also \cite{jacob} Chap. 4).

This example is pushed further by Hiraba \cite{hiraba} who studies the density of the corresponding Markov semi-group using the Malliavin calculus in combining the approaches of \cite{bichteler-gravereaux-jacod} and \cite{leandre1}.
He obtains the existence of a density under the hypothesis that $\alpha(x)$ is $\mathcal{C}^4$ bounded with bounded derivatives and  $0<\lambda_1\leq\alpha(x)\leq\lambda_2<2$ (see also Negoro \cite{negoro} for an analytic proof of this result under $\mathcal{C}^\infty$ assumptions).\\

The method of Dirichlet forms allows several improvements of this subject. First about the correspondence between $\nu(t,x,dy)$ and the pair $(c(t,x,u), N(dt,du))$ it is not necessary to deal with a Poisson measure whose intensity be the Lebesgue measure but only that it carries a Dirichlet form satisfying hypotheses (H0) to (H3). The choice of the operator $\gamma$ is also flexible.

Second the existence of densities may be performed under weaker assumptions on the function $c(t,x,u)$.

Let us write explicitly this example. The operator
$$Af(x)=\int_{\mathbb{R}^d\backslash\{0\}}(f(x+y)-f(x)-D_x f(x).y 1_{|y|<1})\zeta(\alpha(x))\frac{dy}{|y|^{1+\alpha(x)}}$$
represents the symbol (\ref{LL}) if $Ae^{iu.x}=e^{iu.x}(-|u|^{\alpha(x)})$ and this is realized if the function $\zeta$ is such that
$$\zeta(\beta)\int_{\mathbb{R}^d\backslash\{0\}}(1-\cos{\xi.y})\frac{dy}{|y|^{d+\beta}}=|\xi|^\beta$$ what gives (cf \cite{hiraba})
$$\zeta(\beta)=(\sin{\frac{\pi\beta}{2}})\frac{\Gamma(1+\beta)\Gamma(\frac{d+\beta}{2})}{\pi^{\frac{d+1}{2}}\Gamma(\frac{1+\beta}{2})}.$$
As these authors let us except the large jumps for the sake of simplicity. We have to obtain the operator
$$A^0f(x)=\int_{|y|<u_0}(f(x+y)-f(x)-D_xf.y)K(x,dy)$$ with
$$K(x,A)=\zeta(\alpha(x))\int_{S^{d-1}}d\sigma\int_0^\infty1_{A\backslash\{0\}}(r\sigma)\frac{dr}{r^{1+\alpha}}$$ thanks to an SDE of the form
$$X_t=x+\int_0^t\int_{S^{d-1}\times\mathbb{R}_+}c(X_{s-},\sigma,z)\tilde{N}(ds,d\sigma,dz).$$ Here $S^{d-1}$ is the unit sphere in $\mathbb{R}^d$ and $d\sigma$ the area measure.  If we choose the intensity of $N$ to be $dtd\sigma dz$ on $\mathbb{R}_+\times S^{d-1}\times\mathbb{R}_+$ and the function $c(x,\sigma,z)=C(x,z)\sigma$, the condition is that the image of the measure $dz$ on $\mathbb{R}_+$ by the function $C(x,z)$ be the measure $\zeta(\alpha(x))\frac{dr}{r^{1+\alpha(x)}} $ and this yields the function
$$c(x,\sigma,z)=\left(\frac{\alpha(x)z}{\zeta(\alpha(x))}+u_0^{-\alpha(x)}\right)^{-\frac{1}{\alpha(x)}}\;\sigma.$$ Our hypotheses on the bottom structure are fulfilled, we may choose $\gamma[j,j^\ast]=(|x|^2\wedge 1)I$ on $\mathbb{R}^d$. If we suppose
\begin{equation}\label{der}
\alpha\mbox{ of class }\mathcal{C}^1\cap Lip\mbox{ and } 0<\lambda_1\leq\alpha(x)\leq\lambda_2<2
\end{equation}
the assumptions R1a) and R1b) are fulfilled and the condition R1c) has been proved by Hiraba (\cite{hiraba} Remark 3.6 p. 43 {\it et seq.}).

We can conclude that under hypotheses (\ref{der}) there exists a Markov process whose transition semi-group admits $A^0$ as generator and this semi-group possesses a density.
\section{Comments.}

The simplest case of SDE with jumps is the case of L\'evy processes themselves. If $Y$ is a L\'evy process with values in
$\mathbb{R}^d$ the method gives
$$\Gamma[Y_t,Y_t^\ast]=\sum_{\alpha\in JT}\gamma[j,j^\ast](\Delta Y_\alpha)$$ where  $j$ is the identity. In this case
it does not seem that one could do better than the Sato criterion \cite{sato1} (or \cite{sato2} Thm 27.7). This induces the following natural question: when a L\'evy measure $\nu$ is carried by a Lipschitzian curve in $\mathbb{R}^d$ and carries a Dirichlet form satisfying $\sum_{s\leq t}\gamma[j,j^\ast](\Delta Y_s)>0$ a.s. is it necessary such that  $\tilde{\nu}^{\star n}\ll\lambda_d$ for some  $n$ where $\tilde{\nu}=(|x|^2\wedge 1)\nu$ ?

When dealing with SDE's driven by L\'evy processes our approach supposes some regularity for the L\'evy measure of the driver process because of the existence of the bottom Dirichlet form and technical condition (H0). It does not need this L\'evy measure possesses a density (cf \cite{bouleau-denis} \S 2.3 and Example 2 above). Then the existence of density for the solution is obtained under weaker hypotheses than those of  L\'eandre (cf \cite{leandre1} \cite{leandre2}) because there is no growth condition for the L\'evy measure near the origin, and also because we do not need that the  measures on  $\mathbb{R}$  from which the L\'evy measure is a sum of images,  have  $\mathcal{C}^1$-densities. A similar remark may be done when comparing our hypotheses with those of Ishikawa and Kunita \cite{ishikawa-kunita} who suppose non degeneracy of the L\'evy measure of the driver L\'evy process. As said in the introduction, the main advantage of the Dirichlet forms method is to allow only $\mathcal{C}^1\cap Lip$ coefficients.
Now these authors obtain also smoothness results that we do not discuss in the present paper.

In another work we  are studying the extension of our arguments to smoothness results thanks to the fact that the gradient defined by formula (\ref{formulegradient}) may be easily iterated.

L\'evy processes and random Poisson measures do possess strong regularizing properties due to the fact that the jumps are independent of the strict past (cf  \cite{bouleau-denis} examples 5.1 to 5.3). This phenomenon has been deepened by Fournier and Giet  \cite{fournier-giet} who obtained density for the solution of an SDE driven by a L\'evy process supposing only an absolutely continuous L\'evy measure and without using the Malliavin calculus. Even if their hypotheses on the coefficients are slightly stronger than ours, this shows that the Malliavin calculus in the spirit of \cite{bichteler-gravereaux-jacod} which seems to be the most powerful  for this aim, have to be crossed with other techniques in order to capture all regularizing properties of L\'evy processes.

In this perspective, we believe that the simplification brought by the lent particle method gives a tool easier to adapt with various arguments.

 Ecole des Ponts,\\
ParisTech, Paris-Est\\
6 Avenue Blaise Pascal\\ 77455 Marne-La-Vallée Cedex 2
FRANCE\\bouleau@enpc.fr \\  \\ Equipe Analyse et Probabilités,
\\Universit\'{e} d'Evry-Val-d'Essonne,\\Boulevard François Mitterrand\\
91025 EVRY Cedex FRANCE\\ldenis@univ-evry.fr
 \end{document}